\documentclass[12pt]{elsarticle}
\journal{Journal of Differential Equations}

\begin{document}

\begin{frontmatter}

\title{Application of Group Analysis to Classification of Systems of Three Second-Order
Ordinary Differential Equations}


\author[NU]{S. Suksern}

\address[NU]{Department of Mathematics, Faculty of Science,
Naresuan University, Phitsanulok 65000, Thailand}

\ead{supapornsu@nu.ac.th}

\author[DUT]{S. Moyo}

\address[DUT]{Durban University of Technology, Department of Mathematics, Statistics
and Physics \& Institute for Systems Science, P O Box 1334, Steve
Biko Campus, Durban 4000, South Africa}

\ead{moyos@dut.ac.za}

\author[SUT]{S.V. Meleshko}

\address[SUT]{Suranaree University of Technology, School of Mathematics, Nakhon
Ratchasima 30000, Thailand}

\ead{sergey@math.sut.ac.th}

\end{frontmatter}

\begin{abstract}
Here we give a complete group classification of the general case of
linear systems of three second-order ordinary differential equations
excluding the case of systems which are studied in the literature.
This is given as the initial step in the study of nonlinear systems
of three second-order ordinary differential equations. In addition
the complete group classification of a system of three linear second-order
ordinary differential equations is done. Four cases of linear systems
of equations are obtained. \end{abstract}

\section{Introduction}

Systems of ordinary differential equations appear in the modeling
of natural phenomena and have generated sufficient interest in theoretical
studies of these equations and their symmetry properties. For instance,
one of the popular problems in the analysis of differential equations
in the 19th century was that of group classification of ordinary differential
equations (see the works by S. Lie \cite{bk:Lie[1883],bk:Lie[1884]}).
Group classification of differential equations here means to classify
given differential equations with respect to arbitrary elements.

In this paper we consider the complete group classification of systems
of three linear second-order ordinary differential equations. Systems
of second-order ordinary differential equations appear in the study
of many applications. The study of their symmetry structure constitutes
an important field of application of the group symmetry analysis method.
Earlier studies in this area initiated by Lie \cite{bk:Lie[1891b]}
gave a complete group classification of a scalar ordinary differential
equation of the form $y^{\prime\prime}=f(x,y)$. Later on L.V.Ovsiannikov
\cite{bk:Ovsiannikov[2004]} did the group classification using a
different approach. In this approach the classification was obtained
by directly solving the determining equations and exploiting the equivalence
transformations. In more recent works the same approach was used in
\cite{bk:Phauk[2013]} for the group classification of more general
types of equations. We note here that in the general case of a scalar
ordinary differential equation, $y^{\prime\prime}=f(x,y,y^{\prime})$,
the application of the method that involves directly solving the determining
equations gives rise to some difficulties. The group classification
of such equations \cite{bk:MahomedLeach[1989]} is based on the enumeration
of all possible Lie algebras of operators acting on the plane $(x,y)$.
In the works by Lie \cite{bk:Lie[1891b]} is the classification of
all non similar Lie algebras (under complex change of variables) in
two complex variables. In 1992, Gonzalez-Lopez et al. ordered the
Lie classification of realizations of complex Lie algebras and extended
it to the real case \cite{bk:GonzalezKamranOlver[1992a]}. In the
literature and references therein is a large amount of results on
the dimension and structure of symmetry algebras of linearizable ordinary
differential equations (see \cite{bk:MahomedLeach[1989],bk:MahomedLeach[1990],bk:HandbookLie_v3,bk:Gorringeleach[1988],bk:Ovsiannikov[1978],bk:WafoMahomed[2000]}).

It is apparent from these sources that there is a significant number
of studies on symmetry properties of scalar ordinary differential
equations but equally not as much on the group classification of systems
of three linear second-order equations. In recent works \cite{bk:BoykoPopovychShapoval[2012],bk:WafoSoh[2010],bk:Meleshko[2011],bk:Campoamor-Stursberg[2011],bk:Campoamor-Stursberg[2012]}
the authors focused on the study of systems of second-order ordinary
differential equations with constant coefficients of the form
\begin{equation}
\mathbf{y}^{\prime\prime}=M\mathbf{y},\label{feb27.1}
\end{equation}
 where $M$ is a matrix with constant entries.

In the general case of systems of two linear second-order ordinary
differential equations the more advanced results are obtained in \cite{bk:WafoMahomed[2000]},
where the canonical form
\begin{equation}
\mathbf{y}^{\prime\prime}=\left(\begin{array}{cc}
a(x) & b(x)\\
c(x) & -a(x)
\end{array}\right)\mathbf{y},\label{eq:oct09.2013.1}
\end{equation}
 was used by the authors to obtain several admitted Lie groups. We
note that the list of all distinguished representatives of systems
of two linear second-order ordinary differential equations found in
\cite{bk:WafoMahomed[2000]} was not exhaustive and hence this formed
the basis of the paper \cite{bk:MoyoMeleshkoOguis[2013]} where the
complete group classification of two linear second-order ordinary
differential equations using Ovsiannikov's approach was performed.
It is also worth to note here that the form (\ref{eq:oct09.2013.1})
allows one to apply an algebraic approach for group classification.
The algebraic approach takes into account algebraic properties of
an admitted Lie group and the knowledge of the algebraic structure
of admitted Lie algebras that can significantly simplify the group
classification. In particular, the group classification of a single
second-order ordinary differential equation, done by the founder of
the group analysis method, S.Lie \cite{bk:Lie[1883],bk:Lie[1891b]},
cannot be performed without using the algebraic structure of admitted
Lie groups. Recently the algebraic properties for group classification
was applied in \cite{bk:BihloBihloPopovych[2011],bk:BihloBihloPopovych[2012],bk:PopovychBihlo[2012],bk:PopovychIvanovaEshraghi[2004],bk:PopovychKunzingerEshraghi[2004],bk:Chirkunov[2012],bk:Kasatkin[2012],bk:GrigorievMeleshkoSuriyawichitseranee[2013]}.
We also note that the use of the algebraic structure of admitted Lie
groups completely simplified the group classification of equations
describing behavior of fluids with internal inertia in \cite{bk:SiriwatMeleshko[2014],bk:VorakaMeleshko[2014]}.
The study of \cite{bk:MkhizeMoyoMeleshko[2014]} showed that the
problem of classification of systems of two linear second-order ordinary
differential equations using the algebraic approach leads to
the study of more cases than those found in \cite{bk:MoyoMeleshkoOguis[2013]},
where Ovsiannikov's approach was applied.

The system considered in the current paper is a generalization; it
is a system of three linear second-order ordinary differential equations.
We exclude from our consideration the study of systems of second-order
ordinary differential equations with constant coefficients and the
degenerate case given as follows:
\[
y^{\prime\prime}=F(x,y,z,u),\ z^{\prime\prime}=G(x,y,z,u),\ u^{\prime\prime}=0,
\]
 or
\[
y^{\prime\prime}=F(x,y,z,u),\ z^{\prime\prime}=G(x,y,z),\ u^{\prime\prime}=H(x,y,z)=0.
\]

The results found here are new and have not been reported in the literature
as far as we are aware.

The paper is organized as follows:

The first part of the paper deals with the preliminary study of systems
of three second-order nonlinear equations of the form $y^{\prime\prime}=F(x,y,z,u),z^{\prime\prime}=G(x,y,z,u)$
and $u^{\prime\prime}=H(x,y,z,u)$ using Ovsiannikov's approach \cite{bk:Ovsiannikov[2004]}.
This approach involves simplifying one generator and finding the associated
functions. These functions are then used to solve the determining
equations. Using this approach one can classify second-order nonlinear
differential equations but this case has been left for a current project.
The second part of the paper gives a complete treatment of systems
of second-order linear equations. This is then followed by the results
and conclusion.

\section{Preliminary study of systems of nonlinear equations}

A system of three second-order nonlinear differential equations of
the form
\begin{equation}
y^{\prime\prime}=F(x,y,z,u),\ \ z^{\prime\prime}=G(x,y,z,u),\ \ u^{\prime\prime}=H(x,y,z,u)\label{eq:main}
\end{equation}
 is considered in this section. In matrix form equations (\ref{eq:main})
are given by
\begin{equation}
\label{eq:oct5.2013.1}
\mathbf{y}^{\prime\prime}=\mathbf{F(}x,\mathbf{y)},
\end{equation}
 where
\[
\mathbf{y}=\left(\begin{array}{c}
y\\
z\\
u
\end{array}\right),\;\;\mathbf{F}=\left(\begin{array}{c}
F(x,y,z,u)\\
G(x,y,z,u)\\
H(x,y,z,u)
\end{array}\right).
\]

We consider a system of nonlinear equations here as it will later
allow us to separate equations given into their respective classes.
We exclude from the study the degenerate systems which are equivalent
with respect to change of the dependent and independent variables
to one of the classes: (a) the class where $F_{z}=0,\,\,\, F_{u}=0$;
(b) the class where $F_{u}=0,\,\,\, G_{u}=0$. The class (a) is characterized
by the property that one of the equations has the form $y^{\prime\prime}=F(x,y)$.
For the class (b) two equations are of the form
\[
y^{\prime\prime}=F(x,y,z),\ \ z^{\prime\prime}=G(x,y,z).
\]
 Hence for these classes the group classification is reduced to the
study of a single equation or a system of two equations. It is noted
that group classifications of linear equations with the number of
equations $n\leq2$ is complete.

\subsection{Equivalence transformations}

Calculations show that the equivalence Lie group is defined by the
generators:
\[
\begin{array}{c}
X_{1}^{e}=y\partial_{y}+F\partial_{F},\ \ X_{2}^{e}=z\partial_{y}+G\partial_{F},\ \ X_{3}^{e}=u\partial_{y}+H\partial_{F},\\
X_{4}^{e}=y\partial_{z}+F\partial_{G},\ \ X_{5}^{e}=z\partial_{z}+G\partial_{G},\ \ X_{6}^{e}=u\partial_{z}+H\partial_{G},\\
X_{7}^{e}=y\partial_{u}+F\partial_{H},\ \ X_{8}^{e}=z\partial_{u}+G\partial_{H},\ \ X_{9}^{e}=u\partial_{u}+H\partial_{H},\\
X_{10}^{e}=\phi_{1}(x)\partial_{y}+\phi_{1}^{\prime\prime}(x)\partial_{F},\ \ X_{11}^{e}=\phi_{2}(x)\partial_{z}+\phi_{2}^{\prime\prime}(x)\partial_{G},\\
X_{12}^{e}=\phi_{3}(x)\partial_{u}+\phi_{3}^{\prime\prime}(x)\partial_{H},\\
X_{13}^{e}=2\xi(x)\partial_{x}+\xi^{\prime}(x)(y\partial_{y}+z\partial_{z}+u\partial_{u}\\
-3F\partial_{F}-3G\partial_{G}-3H\partial_{H})+\xi^{\prime\prime\prime}(x)(y\partial_{F}+z\partial_{G}+u\partial_{H}),
\end{array}
\]
 where $\xi(x)$ and $\phi_{i}(x)$, $(i=1,2,3)$ are arbitrary functions.

The transformations related with the generators $X_{i}^{e}$, $(i=1,2,...,9)$
correspond to the linear change of the dependent variables $\widetilde{\mathbf{y}}=P\mathbf{y}$
with a constant nonsingular matrix $P$. The transformations corresponding
to the generators $X_{i}^{e}$, $(i=10,11,12)$ define the change
\[
\widetilde{y}=y+\varphi_{1}(x),\quad\widetilde{z}=z+\varphi_{2}(x),\quad\widetilde{u}=u+\varphi_{3}(x).
\]
 The equivalence transformation related with the generator $X_{13}^{e}$
is
\[
\widetilde{x}=\varphi(x),\:\widetilde{y}=y\psi(x),\:\widetilde{z}=z\psi(x),\:\widetilde{u}=u\psi(x),
\]
 where the functions $\varphi(x)$ and $\psi(x)$ satisfy the condition
\begin{equation}
\frac{\varphi^{\prime\prime}}{\varphi^{\prime}}=2\frac{\psi^{\prime}}{\psi}.\label{eq:feb27}
\end{equation}

\subsection{Determining equations}

Consider the generator
\[
X=\xi(x,y,z,u)\frac{\partial}{\partial x}+\eta_{1}(x,y,z,u)\frac{\partial}{\partial y}+\eta_{2}(x,y,z,u)\frac{\partial}{\partial z}+\eta_{3}(x,y,z,u)\frac{\partial}{\partial u}.
\]
 According to the Lie algorithm \cite{bk:Ovsiannikov[1978]}, $X$
is admitted by system (\ref{eq:main}) if it satisfies the associated
determining equations. The first part of the determining equation
is given by
\[
3F\xi_{1}+G\xi_{2}+H\xi_{3}=3\xi_{1}^{\prime\prime}y+\xi_{2}^{\prime\prime}z+\xi_{3}^{\prime\prime}u-\xi_{0}^{\prime\prime}-2h_{1}^{\prime},
\]
\[
F\xi_{1}+3G\xi_{2}+H\xi_{3}=\xi y+3\xi_{2}^{\prime\prime}z+\xi_{3}^{\prime\prime}u-\xi_{0}^{\prime\prime}-2h_{2}^{\prime},
\]
\[
F\xi_{1}+G\xi_{2}+3H\xi_{3}=\xi_{1}^{\prime\prime}y+\xi_{2}^{\prime\prime}z+3\xi_{3}^{\prime\prime}u-\xi_{0}^{\prime\prime}-2h_{3}^{\prime},
\]
\[
F\xi_{2}=\xi_{2}^{\prime\prime}y+h_{12}^{\prime},\,\,\, F\xi_{3}=\xi_{3}^{\prime\prime}y+h_{13}^{\prime},
\]
\[
G\xi_{1}=\xi_{1}^{\prime\prime}z+h_{21}^{\prime},\,\,\, G\xi_{3}=\xi_{3}^{\prime\prime}z+h_{23}^{\prime},
\]
\[
H\xi_{1}=\xi_{1}^{\prime\prime}u+h_{31}^{\prime},\,\,\, H\xi_{2}=\xi_{2}^{\prime\prime}u+h_{32}^{\prime}
\]
 where
\[
\xi(x,y,z,u)=\xi_{1}(x)y+\xi_{2}(x)z+\xi_{3}(x)u+\xi_{0}(x)
\]
 and the functions $h_{i}=h_{i}(x)$, and $h_{ij}=h_{ij}(x)$.

From these equations one can conclude that $\xi_{1}^{2}+\xi_{2}^{2}+\xi_{3}^{2}\neq0$
only for the case where two of the equations are equivalent to the
free particle equation, for instance,
\[
G=0,\ \ \ H=0.
\]

Hence we consider the case where
\[
\xi_{1}=0,\ \ \xi_{2}=0,\ \ \xi_{3}=0.
\]
 The determining equations in this case are given by
\[
\begin{array}{c}
F_{y}(y(\xi^{\prime}+k_{1})+zk_{2}+uk_{3}+\zeta_{1})+F_{z}(yk_{4}+(\xi^{\prime}+k_{5})z+uk_{6}+\zeta_{2})\\
+F_{u}(yk_{7}+zk_{8}+(\xi^{\prime}+k_{9})u+\zeta_{3})\\
+2F_{x}\xi-\xi^{\prime\prime\prime}y+(3\xi^{\prime}-k_{1})F-k_{2}G-k_{3}H-\zeta_{1}^{\prime\prime}=0,
\end{array}
\]
\[
\begin{array}{c}
G_{y}(y(\xi^{\prime}+k_{1})+zk_{2}+uk_{3}+\zeta_{1})+G_{z}(yk_{4}+(\xi^{\prime}+k_{5})z+uk_{6}+\zeta_{2})\\
+G_{u}(yk_{7}+zk_{8}+(\xi^{\prime}+k_{9})u+\zeta_{3})\\
+2G_{x}\xi-\xi^{\prime\prime\prime}z-k_{4}F+(3\xi^{\prime}-k_{5})G-k_{6}H-\zeta_{2}^{\prime\prime}=0,
\end{array}
\]
\[
\begin{array}{c}
H_{y}(y(\xi^{\prime}+k_{1})+zk_{2}+uk_{3}+\zeta_{1})+H_{z}(yk_{4}+(\xi^{\prime}+k_{5})z+uk_{6}+\zeta_{2})\\
+H_{u}(yk_{7}+zk_{8}+(\xi^{\prime}+k_{9})u+\zeta_{3})\\
+2H_{x}\xi-\xi^{\prime\prime\prime}u-k_{7}F-k_{8}G+(3\xi^{\prime}-k_{9})H-\zeta_{3}^{\prime\prime}=0,
\end{array}
\]
 where an admitted generator has the form
\[
\begin{array}{c}
X=2\xi\frac{\partial}{\partial x}+((\xi^{\prime}+k_{1})y+k_{2}z+k_{3}u+\zeta_{1})\frac{\partial}{\partial y}\\
+(k_{4}y+(\xi^{\prime}+k_{5})z+k_{6}u+\zeta_{2})\frac{\partial}{\partial z}\\
+(k_{7}y+k_{8}z+(\xi^{\prime}+k_{9})u+\zeta_{3})\frac{\partial}{\partial u}
\end{array}
\]
 with $\xi=\xi(x)$, $\zeta_{i}=\zeta_{i}(x)$ and $k_{i}$, $i=1,2,3$
constant.

For further analysis the study of the determining equations is separated
into the cases:

(a) the case in which there is at least one admitted generator with
$\xi\neq0$;

(b) the case in which for all admitted generators $\xi=0$.

\subsection{Case $\xi\neq0$}

We consider the generator $X_{o}$ for which $\xi\neq0$. Using the
equivalence transformation:
\[
y_{1}=y+\varphi_{1}(x),\ \ z_{1}=z+\varphi_{2},\ \ u_{1}=u+\varphi_{3}(x),
\]
 the generator $X_{o}$ becomes
\[
\begin{array}{c}
X_{o}=2\xi\frac{\partial}{\partial x}+((\xi^{\prime}+k_{1})y_{1}+k_{2}z_{1}+k_{3}u_{1}+\widetilde{\zeta}_{1})\frac{\partial}{\partial y_{1}}\\
+(k_{4}y_{1}+(\xi^{\prime}+k_{5})z_{1}+k_{6}u_{1}+\widetilde{\zeta}_{2})\frac{\partial}{\partial z_{1}}\\
+(k_{7}y_{1}+k_{8}z_{1}+(\xi^{\prime}+k_{9})u_{1}+\widetilde{\zeta}_{3})\frac{\partial}{\partial u_{1}},
\end{array}
\]
 where
\[
\widetilde{\zeta}_{1}=2\xi\varphi_{1}^{\prime}-(\xi^{\prime}+k_{1})\varphi_{1}-k_{2}\varphi_{2}-k_{3}\varphi_{3}+\zeta_{1},
\]
\[
\widetilde{\zeta}_{2}=2\xi\varphi_{2}^{\prime}-k_{4}\varphi_{1}-(\xi^{\prime}+k_{5})\varphi_{2}-k_{6}\varphi_{3}+\zeta_{2},
\]
\[
\widetilde{\zeta}_{3}=2\xi\varphi_{3}^{\prime}-k_{7}\varphi_{1}-k_{8}\varphi_{2}-(\xi^{\prime}+k_{9})\varphi_{3}+\zeta_{3}.
\]
 The functions $\varphi_{i}(x)$ ($i=1,2,3$) can be chosen such that
\[
\widetilde{\zeta}_{1}=0,\,\,\,\widetilde{\zeta}_{2}=0,\,\,\,\widetilde{\zeta}_{3}=0.
\]
 Then without loss of generality one can assume that the generator
$X_{o}$ has the form
\[
\begin{array}{c}
X_{o}=2\xi\frac{\partial}{\partial x}+((\xi^{\prime}+k_{1})y+k_{2}z+k_{3}u)\frac{\partial}{\partial y}+(k_{4}y+(\xi^{\prime}+k_{5})z+k_{6}u)\frac{\partial}{\partial z}\\
+(k_{7}y+k_{8}z+(\xi^{\prime}+k_{9})u)\frac{\partial}{\partial u}.
\end{array}
\]

The equivalence transformation
\[
x_{1}=\alpha(x),\ \ y_{1}=y\beta(x),\ \ z_{1}=z\beta(x),\ \ u_{1}=u\beta(x),
\]
 where
\[
\alpha^{\prime\prime}\beta=2\alpha^{\prime}\beta^{\prime},\;\;\;(\alpha^{\prime}\beta\neq0),
\]
 reduces the generator $X_{o}$ to
\[
\begin{array}{c}
X_{o}=2\alpha^{\prime}\xi\frac{\partial}{\partial x_{1}}+((2\xi\beta^{\prime}/\beta+\xi^{\prime}+k_{1})y_{1}+k_{2}z_{1}+k_{3}u_{1})\frac{\partial}{\partial y_{1}}\\
+(k_{4}y_{1}+(2\xi\beta^{\prime}/\beta+\xi^{\prime}+k_{5})z_{1}+k_{6}u_{1})\frac{\partial}{\partial z_{1}}\\
+(k_{7}y_{1}+k_{8}z_{1}+(2\xi\beta^{\prime}/\beta+\xi^{\prime}+k_{9})u_{1})\frac{\partial}{\partial u_{1}}.
\end{array}
\]
 Choosing $\beta(x)$ such that $2\xi\beta^{\prime}/\beta+\xi^{\prime}=0$
makes the generator $X_{o}$ become
\[
\begin{array}{c}
X_{o}=2\alpha^{\prime}\xi\frac{\partial}{\partial x_{1}}+(k_{1}y_{1}+k_{2}z_{1}+k_{3}u_{1})\frac{\partial}{\partial y_{1}}+(k_{4}y_{1}+k_{5}z_{1}+k_{6}u_{1})\frac{\partial}{\partial z_{1}}\\
+(k_{7}y_{1}+k_{8}z_{1}+k_{9}u_{1})\frac{\partial}{\partial u_{1}}.
\end{array}
\]
 Note that in this case
\[
\frac{d(\alpha^{\prime}\xi)}{dx_{1}}=0.
\]
 This means
\[
\frac{d(\alpha^{\prime}\xi)}{dx_{1}}=\frac{(\alpha^{\prime}\xi)^{\prime}}{\alpha^{\prime}}=\xi^{\prime}+\frac{\alpha^{\prime\prime}}{\alpha^{\prime}}\xi=-2\xi\frac{\beta^{\prime}}{\beta}+2\frac{\beta^{\prime}}{\beta}\xi=0.
\]
 Hence without loss of generality we can assume that the generator
$X_{o}$ has the form
\[
\begin{array}{c}
X_{o}=k\frac{\partial}{\partial x}+(k_{1}y+k_{2}z+k_{3}u)\frac{\partial}{\partial y}+(k_{4}y+k_{5}z+k_{6}u)\frac{\partial}{\partial z}\\
+(k_{7}y+k_{8}z+k_{9}u)\frac{\partial}{\partial u},
\end{array}
\]
 where $k=2\alpha^{\prime}\xi\neq0$ is constant. We rewrite the generator
$X_{o}$ in the form
\[
\begin{array}{c}
X_{o}=\partial_{x}+(a_{11}y+a_{12}z+a_{13}u)\frac{\partial}{\partial y}+(a_{21}y+a_{22}z+a_{23}u)\frac{\partial}{\partial z}\\
+(a_{31}y+a_{32}z+a_{33}u)\frac{\partial}{\partial u}.
\end{array}
\]
 The determining equations become
\begin{equation}
\begin{array}{l}
\left(a_{11}y+a_{12}z+a_{13}u\right)F_{y}+\left(a_{21}y+a_{22}z+a_{23}u\right)F_{z}\\
+\left(a_{31}y+a_{32}z+a_{33}u\right)F_{u}+F_{x}=a_{11}F+a_{12}G+a_{13}H,\\
\left(a_{11}y+a_{12}z+a_{13}u\right)G_{y}+\left(a_{21}y+a_{22}z+a_{23}u\right)G_{z}\\
+\left(a_{31}y+a_{32}z+a_{33}u\right)G_{u}+G_{x}=a_{21}F+a_{22}G+a_{23}H,\\
\left(a_{11}y+a_{12}z+a_{13}u\right)H_{y}+\left(a_{21}y+a_{22}z+a_{23}u\right)H_{z}\\
+\left(a_{31}y+a_{32}z+a_{33}u\right)H_{u}+H_{x}=a_{31}F+a_{32}G+a_{33}H.
\end{array}\label{cl:01}
\end{equation}
 Here $a_{ij}$, ($i,j=1,2,3$) are constant. In matrix form these
equations are rewritten as
\begin{equation}
\mathbf{F}_{x}+\left((A\mathbf{y})\cdot\nabla\right)\mathbf{F}-A\mathbf{F}=0,\label{eq:aug1011.1}
\end{equation}
 where
\[
A=\left(\begin{array}{ccc}
a_{11} & a_{12} & a_{13}\\
a_{21} & a_{22} & a_{23}\\
a_{31} & a_{32} & a_{33}
\end{array}\right),\;\nabla=\left(\begin{array}{c}
\partial_{y}\\
\partial_{z}\\
\partial_{u}
\end{array}\right).
\]
 Here {}``$\cdot$'' means the scalar product
\[
\mathbf{b}\cdot\nabla=b_{i}\partial_{y_{i}},
\]
 where a vector $\mathbf{b}=(b_{1},b_{2},b_{3}),$ $\mathbf{y}=(y_{1},y_{2},y_{3})$,
and it is also used standard agreement: summation with respect to
a repeat index.

Further simplifications are related with simplifications of the matrix
$A$.

We apply the change $\widetilde{\mathbf{y}}=P\mathbf{y}$ where $P$
is a nonsingular matrix with constant entries
\[
P=\left(\begin{array}{ccc}
p_{11} & p_{12} & p_{13}\\
p_{21} & p_{22} & p_{23}\\
p_{31} & p_{32} & p_{33}
\end{array}\right).
\]
 Equations (\ref{eq:main}) become
\[
\widetilde{\mathbf{y}}^{\prime\prime}=\widetilde{\mathbf{F}}(x,\widetilde{\mathbf{y}})
\]
 where
\[
\widetilde{\mathbf{F}}(x,\widetilde{\mathbf{y}})=P\mathbf{F}(x,P^{-1}\widetilde{\mathbf{y}}).
\]
 The partial derivatives $\partial_{y},\partial_{z}$ and $\partial_{u}$
are changed as follows:
\[
\nabla=P^{t}\widetilde{\nabla}.
\]
 Hence equations (\ref{eq:aug1011.1}) become
\[
\begin{array}{l}
\left((AP^{-1}\widetilde{\mathbf{y}})\cdot P^{t}\widetilde{\nabla}\right)(P^{-1}\widetilde{\mathbf{F}})+P^{-1}\widetilde{\mathbf{F}}_{x}-AP^{-1}\widetilde{\mathbf{F}}\\
=P^{-1}\left(\left((PAP^{-1}\widetilde{\mathbf{y}})\cdot\widetilde{\nabla}\right)\widetilde{\mathbf{F}}+\widetilde{\mathbf{F}}_{x}-PAP^{-1}\widetilde{\mathbf{F}}\right)\\
=P^{-1}\left(\left((\widetilde{A}\widetilde{\mathbf{y}})^{t}\widetilde{\nabla}\right)\cdot\widetilde{\mathbf{F}}+\widetilde{\mathbf{F}}_{x}-\widetilde{A}\widetilde{\mathbf{F}}\right)=0
\end{array}
\]
 where
\[
\widetilde{A}=PAP^{-1}.
\]
 This means that the change $\widetilde{\mathbf{y}}=P\mathbf{y}$
reduces equation (\ref{eq:aug1011.1}) to the same form with the matrix
$A$ changed. The infinitesimal generator is also changed as
\[
X_{o}=\partial_{x}+(\widetilde{A}\widetilde{\mathbf{y}})\cdot\widetilde{\nabla}.
\]
 Using this change matrix $A$ can be presented in the Jordan form.
For a real-valued $3\times3$ matrix $A$, if the matrix $P$ also
has real-valued entries, then the Jordan matrix is one of the following
four types:
\begin{equation}
J_{1}=\left(\begin{array}{ccc}
a & 0 & 0\\
0 & b & 0\\
0 & 0 & d
\end{array}\right),\;\; J_{2}=\left(\begin{array}{ccc}
a & 0 & 0\\
0 & b & c\\
0 & -c & b
\end{array}\right),\;\; J_{3}=\left(\begin{array}{ccc}
a & 0 & 0\\
0 & b & 1\\
0 & 0 & b
\end{array}\right),\;\; J_{4}=\left(\begin{array}{ccc}
a & 1 & 0\\
0 & a & 1\\
0 & 0 & a
\end{array}\right),\label{eq:Jordan_form}
\end{equation}
 where $a$, $b$, $c$ and $d>0$ are real numbers.

\subsubsection{Case $A=J_{1}$}

We assume that
\[
A=\left(\begin{array}{ccc}
a & 0 & 0\\
0 & b & 0\\
0 & 0 & d
\end{array}\right).
\]
 In this case the equations for the functions F, G and H are
\[
\begin{array}{l}
ayF_{y}+bzF_{z}+duF_{u}+F_{x}=aF,\\
ayG_{y}+bzG_{z}+duG_{u}+G_{x}=bG,\\
ayH_{y}+bzH_{z}+duH_{u}+H_{x}=dH.
\end{array}
\]
 The general solution of these equations is
\begin{equation}
\begin{array}{c}
F\left(x,y,z,u\right)=e^{ax}f\left(s,v,w\right),\,\, G\left(x,y,z,u\right)=e^{bx}g\left(s,v,w\right),\\
H\left(x,y,z,u\right)=e^{dx}h\left(s,v,w\right)
\end{array}\label{cl:02}
\end{equation}
 where
\[
s=ye^{-ax},\,\, v=ze^{-bx},\,\, w=ue^{-dx}.
\]
 The admitted generator is
\[
X_{0}=\frac{\partial}{\partial x}+ay\frac{\partial}{\partial y}+bz\frac{\partial}{\partial z}+du\frac{\partial}{\partial u}.
\]


\subsubsection{Case $A=J_{2}$}

We assume that
\[
A=\left(\begin{array}{ccc}
a & 0 & 0\\
0 & b & c\\
0 & -c & b
\end{array}\right).
\]
 In this case equations (\ref{cl:01}) become
\begin{equation}
\begin{array}{l}
ayF_{y}+\left(bz+cu\right)F_{z}+\left(-cz+bu\right)F_{u}+F_{x}=aF,\\
ayG_{y}+\left(bz+cu\right)G_{z}+\left(-cz+bu\right)G_{u}+G_{x}=bG+cH,\\
ayH_{y}+\left(bz+cu\right)H_{z}+\left(-cz+bu\right)H_{u}+H_{x}=cG+bH.
\end{array}\label{cl:03}
\end{equation}
 Here again we introduce the variables
\[
\begin{array}{c}
s=ye^{-ax},\,\, v=e^{-bx}\left(z\cos\left(cx\right)-u\sin\left(cx\right)\right),\\
w=e^{-bx}\left(z\sin\left(cx\right)+u\cos\left(cx\right)\right),
\end{array}
\]
 equations (\ref{cl:03}) become
\[
F_{x}-aF=0,\,\, G_{x}-bG-cH=0,\,\, H_{x}+cG-bH=0.
\]
 The general solution of these equations is
\begin{equation}
\begin{array}{l}
F\left(x,y,z,u\right)=e^{ax}f\left(s,v,w\right),\,\,\\
G\left(x,y,z,u\right)=e^{bx}\left(\cos\left(cx\right)g\left(s,v,w\right)+\sin\left(cx\right)h\left(s,v,w\right)\right),\,\,\\
H\left(x,y,z,u\right)=e^{bx}\left(-\sin\left(cx\right)g\left(s,v,w\right)+\cos\left(cx\right)h\left(s,v,w\right)\right)
\end{array}\label{cl:04}
\end{equation}
 where $f(s,v,w),g(s,v,w)$ and $h(s,v,w)$ are arbitrary functions.
The admitted generator is
\[
X_{0}=\frac{\partial}{\partial x}+ay\frac{\partial}{\partial y}+\left(bz+cu\right)\frac{\partial}{\partial z}+\left(-cz+bu\right)\frac{\partial}{\partial u}.
\]


\subsubsection{Case $A=J_{3}$}

We assume that
\[
A=\left(\begin{array}{ccc}
a & 0 & 0\\
0 & b & 1\\
0 & 0 & b
\end{array}\right).
\]
 In this case equations (\ref{cl:01}) become
\begin{equation}
\begin{array}{l}
ayF_{y}+\left(bz+cu\right)F_{z}+buF_{u}+F_{x}=aF,\\
ayG_{y}+\left(bz+cu\right)G_{z}+buG_{u}+G_{x}=bG+H,\\
ayH_{y}+\left(bz+cu\right)H_{z}+buH_{u}+H_{x}=bH.
\end{array}\label{cl:05}
\end{equation}
 As in the previous case we introduce the variables

\[
s=ye^{-ax},\,\, v=e^{-bx}\left(z-ux\right),\,\, w=e^{-bx}u
\]
 so that equations (\ref{cl:05}) become
\[
F_{x}-aF=0,\,\, G_{x}-bG-H=0,\,\, H_{x}-bH=0.
\]
 The general solution of these equations is
\begin{equation}
\begin{array}{l}
F\left(x,y,z,u\right)=e^{ax}f\left(s,v,w\right),\,\,\\
G\left(x,y,z,u\right)=e^{bx}\left(h\left(s,v,w\right)x+g\left(s,v,w\right)\right),\,\,\\
H\left(x,y,z,u\right)=e^{bx}h\left(s,v,w\right),
\end{array}\label{cl:06}
\end{equation}
 where $f(s,v,w),g(s,v,w)$ and $h(s,v,w)$ are arbitrary functions.
The admitted generator is
\[
X_{0}=\frac{\partial}{\partial x}+ay\frac{\partial}{\partial y}+\left(bz+cu\right)\frac{\partial}{\partial z}+bu\frac{\partial}{\partial u}.
\]


\subsubsection{Case $A=J_{4}$}

We assume that
\[
A=\left(\begin{array}{ccc}
a & 1 & 0\\
0 & a & 1\\
0 & 0 & a
\end{array}\right).
\]
 In this case equations (\ref{cl:01}) become
\begin{equation}
\begin{array}{l}
\left(ay+z\right)F_{y}+\left(az+u\right)F_{z}+auF_{u}+F_{x}=aF,\\
\left(ay+z\right)G_{y}+\left(az+u\right)G_{z}+auG_{u}+G_{x}=aG+H,\\
\left(ay+z\right)H_{y}+\left(az+u\right)H_{z}+auH_{u}+H_{x}=aH.
\end{array}\label{cl:07}
\end{equation}
 Introducing the variables
\[
s=e^{-ax}\left(y-xz+\frac{1}{2}x^{2}u\right),\,\, v=e^{-ax}\left(z-xu\right),\,\, w=e^{-ax}u
\]
 so that equations (\ref{cl:07}) become
\[
F_{x}-aF-G=0,\,\, G_{x}-aG-H=0,\,\, H_{x}-aH=0.
\]
 The general solution of these equations is
\begin{equation}
\begin{array}{l}
F\left(x,y,z,u\right)=e^{ax}\left(h\left(s,v,w\right)\frac{x^{2}}{2}+g\left(s,v,w\right)+f\left(s,v,w\right)\right),\,\,\\
G\left(x,y,z,u\right)=e^{ax}\left(h\left(s,v,w\right)x+g\left(s,v,w\right)\right),\,\,\\
H\left(x,y,z,u\right)=e^{bx}h\left(s,v,w\right)
\end{array}\label{cl:08}
\end{equation}
 where $f(s,v,w),g(s,v,w)$ and $h(s,v,w)$ are arbitrary functions.
The admitted generator is
\[
X_{0}=\frac{\partial}{\partial x}+(ay+z)\frac{\partial}{\partial y}+\left(az+u\right)\frac{\partial}{\partial z}+au\frac{\partial}{\partial u}.
\]

\subsection{Case $\xi=0$}

Substituting $\xi=0$ into the determining equations we find that
\begin{equation}
\begin{array}{c}
(a_{11}y+a_{12}z+a_{13}u+\zeta_{1})F_{y}+(a_{21}y+a_{22}z+a_{23}u+\zeta_{2})F_{z}\\
+(a_{31}y+a_{32}z+a_{33}u+\zeta_{3})F_{u}=a_{11}F+a_{12}G+a_{13}H+\zeta_{1}^{\prime\prime},\\
(a_{11}y+a_{12}z+a_{13}u+\zeta_{1})G_{y}+(a_{21}y+a_{22}z+a_{23}u+\zeta_{2})G_{z}\\
+(a_{31}y+a_{32}z+a_{33}u+\zeta_{3})G_{u}=a_{21}F+a_{22}G+a_{23}H+\zeta_{2}^{\prime\prime},\\
(a_{11}y+a_{12}z+a_{13}u+\zeta_{1})H_{y}+(a_{21}y+a_{22}z+a_{23}u+\zeta_{2})H_{z}\\
+(a_{31}y+a_{32}z+a_{33}u+\zeta_{3})H_{u}=a_{31}F+a_{32}G+a_{33}H+\zeta_{3}^{\prime\prime},
\end{array}\label{eq:oct5.5}
\end{equation}
 or in matrix form
\begin{equation}
\left((A\mathbf{y}+\mathbf{h})\cdot\nabla\right)\mathbf{F}=A\mathbf{F}+\mathbf{h}^{\prime\prime}\label{eq:oct2_2013.1}
\end{equation}
 where
\[
A=\left(\begin{array}{ccc}
k_{1} & k_{2} & k_{3}\\
k_{4} & k_{5} & k_{6}\\
k_{7} & k_{8} & k_{9}
\end{array}\right),\ \ \mathbf{h}(x)=\left(\begin{array}{c}
\zeta_{1}(x)\\
\zeta_{2}(x)\\
\zeta_{3}(x)
\end{array}\right).
\]

Similarly to the case where $\xi\neq0$ we use the Jordan forms (\ref{eq:Jordan_form})
of the matrix $A$. The admitted generator has the form
\[
\begin{array}{c}
X_{o}=(k_{1}y+k_{2}z+k_{3}u+\zeta_{1}(x))\partial_{y}+(k_{4}y+k_{5}z+k_{6}u+\zeta_{2}(x))\partial_{z}\\
+(k_{7}y+k_{8}z+k_{9}u+\zeta_{3}(x))\partial_{u}.
\end{array}
\]

\subsubsection{Case $A=J_{1}$}

Assuming that $A=J_{1}$, equations (\ref{eq:oct5.5}) for the functions
$F$ and $G$ are
\[
\begin{array}{c}
(ay+h_{1})F_{y}+(bz+h_{2})F_{z}+(du+h_{3})F_{u}=aF+h_{1}^{\prime\prime},\\
(ay+h_{1})G_{y}+(bz+h_{2})G_{z}+(du+h_{3})G_{u}=bG+h_{2}^{\prime\prime},\\
(ay+h_{1})H_{y}+(bz+h_{2})H_{z}+(du+h_{3})H_{u}=dH+h_{3}^{\prime\prime}.
\end{array}
\]
 The general solution of these equations depends on a value of $a,b$
and $d$:
\begin{itemize}
\item {Case: $a\neq0,\, b\neq0,\, d\neq0$}
\[
\begin{array}{c}
aF+h_{1}^{\prime\prime}=(ay+h_{1})f(x,s,v),\ \ bG+h_{2}^{\prime\prime}=(bz+h_{2})g(x,s,v),\\
dH+h_{3}^{\prime\prime}=(du+h_{3})h(x,s,v),\ \ s=(bz+h_{2})^{a}(ay+h_{1})^{-b},\\
v=(du+h_{3})^{b}(bz+h_{2})^{-d}.
\end{array}
\]

\item {Case: $a\neq0,\, b\neq0,\, d=0$}
\[
\begin{array}{c}
aF+h_{1}^{\prime\prime}=(ay+h_{1})f(x,s,v),\ \ bG+h_{2}^{\prime\prime}=(bz+h_{2})g(x,s,v),\\
H=\frac{h_{3}^{\prime\prime}}{a}\ln(ay+h_{1})+h(x,s,v),\ \ s=(bz+h_{2})^{a}(ay+h_{1})^{-b},\\
v=u-\frac{h_{3}}{b}\ln(bz+h_{2}).
\end{array}
\]

\item {Case: $a\neq0,\, b=0,\, d=0$}
\[
\begin{array}{c}
aF+h_{1}^{\prime\prime}=(ay+h_{1})f(x,s,v),\ \ G=\frac{h_{2}^{\prime\prime}}{a}\ln(ay+h_{1})+g(x,s,v),\\
H=\frac{h_{3}^{\prime\prime}}{a}\ln(ay+h_{1})+h(x,s,v),\ \ s=z-\frac{h_{2}}{a}\ln(ay+h_{1}),\\
v=u-\frac{h_{3}}{a}\ln(ay+h_{1}).
\end{array}
\]

\item {Case: $a=0,\, b=0,\, d=0,h_{1}\neq0$}
\[
\begin{array}{c}
F=\frac{h_{1}^{\prime\prime}}{h_{1}}y+f(x,s,v),\ \ G=\frac{h_{2}^{\prime\prime}}{h_{1}}y+g(x,s,v),\\
H=\frac{h_{3}^{\prime\prime}}{h_{1}}y+h(x,s,v),\ \ s=z-\frac{h_{2}}{h_{1}}y,\ \ v=u-\frac{h_{3}}{h_{1}}y.
\end{array}
\]

\item {Case: $a=0,\, b=0,\, d=0,h_{1}=0,h_{2}\neq0$}
\[
\begin{array}{c}
F=f(x,s),\ \ G=\frac{h_{2}^{\prime\prime}}{h_{2}}z+g(x,s),\ \ H=\frac{h_{3}^{\prime\prime}}{h_{2}}z+h(x,s),\ \ s=u-\frac{h_{3}}{h_{2}}z.\end{array}
\]

\end{itemize}
Here $f(x,s,v),g(x,s,v)$ and $h(x,s,v)$ are arbitrary functions.

\subsection{Case $A=J_{2}$}

In this case equations (\ref{eq:oct5.5}) become
\begin{equation}
\begin{array}{l}
(ay+h_{1})F_{y}+(bz+cu+h_{2})F_{z}+(-cz+bu+h_{3})F_{u}=aF+h_{1}^{\prime\prime},\\
(ay+h_{1})G_{y}+(bz+cu+h_{2})G_{z}+(-cz+bu+h_{3})G_{u}=bG+cH+h_{2}^{\prime\prime},\\
(ay+h_{1})H_{y}+(bz+cu+h_{2})H_{z}+(-cz+bu+h_{3})H_{u}=-cG+bH+h_{3}^{\prime\prime}.
\end{array}\label{eq:oct5.1}
\end{equation}

\begin{itemize}
\item Case: $a\neq0$

Introducing the variables
\[
\begin{array}{c}
y=\overline{y}-\frac{h_{1}}{a},\ \ z=\overline{z}-(b^{2}+c^{2})^{-1}(bh_{2}-ch_{3}),\\
u=\overline{u}-(b^{2}+c^{2})^{-1}(ch_{2}+bh_{3}),
\end{array}
\]
\[
\begin{array}{c}
F=\overline{F}-\frac{h_{1}^{\prime\prime}}{a},\ \ G=\overline{G}-(b^{2}+c^{2})^{-1}(bh_{2}^{\prime\prime}-ch_{3}^{\prime\prime}),\\
H=\overline{H}-(b^{2}+c^{2})^{-1}(ch_{2}^{\prime\prime}+bh_{3}^{\prime\prime})
\end{array}
\]
 equations (\ref{eq:oct5.2013.1}) become
\[
\begin{array}{l}
a\overline{y}\overline{F}_{\overline{y}}+(b\overline{z}+c\overline{u})\overline{F}_{\overline{z}}+(-c\overline{z}+b\overline{u})\overline{F}_{\overline{u}}=a\overline{F},\\
a\overline{y}\overline{G}_{\overline{y}}+(b\overline{z}+c\overline{u})\overline{G}_{\overline{z}}+(-c\overline{z}+b\overline{u})\overline{G}_{\overline{u}}=b\overline{G}+c\overline{H},\\
a\overline{y}\overline{H}_{\overline{y}}+(b\overline{z}+c\overline{u})\overline{H}_{\overline{z}}+(-c\overline{z}+b\overline{u})\overline{H}_{\overline{u}}=-c\overline{G}+b\overline{H}.
\end{array}
\]
 In the variables
\[
\overline{y}=ve^{as},\ \ \ \overline{z}=we^{bs}\sin(cs),\ \ \ \overline{u}=we^{bs}\cos(cs)
\]
 these equations are
\[
\overline{F}_{s}=a\overline{F},\quad\overline{G}_{s}=b\overline{G}+c\overline{H},\quad\overline{H}_{s}=-c\overline{G}+b\overline{H}.
\]
 The general solution of the last set of equations is
\[
\overline{F}(x,y,z,u)=e^{as}f(x,v,w),
\]
\[
\overline{G}(x,y,z,u)=e^{as}\left(\cos(cs)g(x,v,w)+\sin(cs)h(x,v,w)\right),
\]
\[
\overline{H}(x,y,z,u)=e^{bs}\left(-\sin(cs)g(x,v,w)+\cos(cs)h(x,v,w)\right)
\]
 where $f(x,v,w),g(x,v,w)$ and $h(x,v,w)$ are arbitrary functions.

\end{itemize}

\begin{itemize}
\item Case: $a=0$

In this case equations (\ref{eq:oct5.2013.1}) become
\begin{equation}
\begin{array}{l}
h_{1}F_{y}+(bz+cu+h_{2})F_{z}+(-cz+bu+h_{3})F_{u}=h_{1}^{\prime\prime},\\
h_{1}G_{y}+(bz+cu+h_{2})G_{z}+(-cz+bu+h_{3})G_{u}=bG+cH+h_{2}^{\prime\prime},\\
h_{1}H_{y}+(bz+cu+h_{2})H_{z}+(-cz+bu+h_{3})H_{u}=-cG+bH+h_{3}^{\prime\prime}.
\end{array}\label{cl:09}
\end{equation}
\begin{itemize}
\item If $h_{1}\neq0$, introducing the variables
\[
\begin{array}{c}
y=\overline{y}h_{1},\ \ z=\overline{z}-(b^{2}+c^{2})^{-1}(bh_{2}-ch_{3}),\\
u=\overline{u}-(b^{2}+c^{2})^{-1}(ch_{2}+bh_{3}),
\end{array}
\]
\[
\begin{array}{c}
F=\overline{F},\ \ G=\overline{G}-(b^{2}+c^{2})^{-1}(bh_{2}^{\prime\prime}-ch_{3}^{\prime\prime}),\\
H=\overline{H}-(b^{2}+c^{2})^{-1}(ch_{2}^{\prime\prime}+bh_{3}^{\prime\prime})
\end{array}
\]
 equations (\ref{cl:09}) become
\[
\begin{array}{l}
\overline{F}_{\overline{y}}+(b\overline{z}+c\overline{u})\overline{F}_{\overline{z}}+(-c\overline{z}+b\overline{u})\overline{F}_{\overline{u}}=h_{1}^{\prime\prime},\\
\overline{G}_{\overline{y}}+(b\overline{z}+c\overline{u})\overline{G}_{\overline{z}}+(-c\overline{z}+b\overline{u})\overline{G}_{\overline{u}}=b\overline{G}+c\overline{H},\\
\overline{H}_{\overline{y}}+(b\overline{z}+c\overline{u})\overline{H}_{\overline{z}}+(-c\overline{z}+b\overline{u})\overline{H}_{\overline{u}}=-c\overline{G}+b\overline{H}.
\end{array}
\]
 In the variables
\[
\begin{array}{c}
s=\overline{y},\ \ \ \ v=e^{-b\overline{y}}\left(\overline{z}\cos\left(c\overline{y}\right)-\overline{u}\sin\left(c\overline{y}\right)\right),\ \\
w=e^{-b\overline{y}}\left(\overline{z}\sin\left(c\overline{y}\right)+\overline{u}\cos\left(c\overline{y}\right)\right)
\end{array}
\]
 these equations become
\[
\overline{F}_{s}=h_{1}^{\prime\prime},\quad\overline{G}_{s}=b\overline{G}+c\overline{H},\quad\overline{H}_{s}=-c\overline{G}+b\overline{H}.
\]
 The general solution of the last set of equations is
\[
\overline{F}(x,y,z,u)=h_{1}^{\prime\prime}s+f(x,v,w),
\]
\[
\overline{G}(x,y,z,u)=e^{bs}\left(\cos(cs)g(x,v,w)+\sin(cs)h(x,v,w)\right),
\]
\[
\overline{H}(x,y,z,u)=e^{bs}\left(-\sin(cs)g(x,v,w)+\cos(cs)h(x,v,w)\right).
\]

\end{itemize}

\begin{itemize}
\item If $h_{1}=0$, in this case equations (\ref{cl:09}) become
\begin{equation}
\begin{array}{l}
(bz+cu+h_{2})F_{z}+(-cz+bu+h_{3})F_{u}=0,\\
(bz+cu+h_{2})G_{z}+(-cz+bu+h_{3})G_{u}=bG+cH+h_{2}^{\prime\prime},\\
(bz+cu+h_{2})H_{z}+(-cz+bu+h_{3})H_{u}=-cG+bH+h_{3}^{\prime\prime}.
\end{array}\label{cl:10}
\end{equation}
 Introducing the variables
\[
z=\overline{z}-(b^{2}+c^{2})^{-1}(bh_{2}-ch_{3}),\ \ u=\overline{u}-(b^{2}+c^{2})^{-1}(ch_{2}+bh_{3}),
\]
\[
\begin{array}{c}
F=\overline{F},\ \ G=\overline{G}-(b^{2}+c^{2})^{-1}(bh_{2}^{\prime\prime}-ch_{3}^{\prime\prime}),\\
H=\overline{H}-(b^{2}+c^{2})^{-1}(ch_{2}^{\prime\prime}+bh_{3}^{\prime\prime})
\end{array}
\]
 equations (\ref{cl:10}) become
\[
\begin{array}{l}
(b\overline{z}+c\overline{u})\overline{F}_{\overline{z}}+(-c\overline{z}+b\overline{u})\overline{F}_{\overline{u}}=0,\\
(b\overline{z}+c\overline{u})\overline{G}_{\overline{z}}+(-c\overline{z}+b\overline{u})\overline{G}_{\overline{u}}=b\overline{G}+c\overline{H},\\
(b\overline{z}+c\overline{u})\overline{H}_{\overline{z}}+(-c\overline{z}+b\overline{u})\overline{H}_{\overline{u}}=-c\overline{G}+b\overline{H}.
\end{array}
\]
 In the variables
\[
\overline{y}=w,\ \ \ \overline{z}=ve^{bs}\sin(cs),\ \ \ \overline{u}=ve^{bs}\cos(cs)
\]
 these equations become
\[
\overline{F}_{s}=0,\quad\overline{G}_{s}=b\overline{G}+c\overline{H},\quad\overline{H}_{s}=-c\overline{G}+b\overline{H}
\]
 and their general solution is
\[
\overline{F}(x,y,z,u)=f(x,v,w),
\]
\[
\overline{G}(x,y,z,u)=e^{bs}\left(\cos(cs)g(x,v,w)+\sin(cs)h(x,v,w)\right),
\]
\[
\overline{H}(x,y,z,u)=e^{bs}\left(-\sin(cs)g(x,v,w)+\cos(cs)h(x,v,w)\right).
\]

\end{itemize}
\end{itemize}

\subsubsection{Case $A=J_{3}$}

In this case equations (\ref{eq:oct5.5}) become
\begin{equation}
\begin{array}{l}
(ay+h_{1})F_{y}+(bz+u+h_{2})F_{z}+(bu+h_{3})F_{u}=aF+h_{1}^{\prime\prime},\\
(ay+h_{1})G_{y}+(bz+u+h_{2})G_{z}+(bu+h_{3})G_{u}=bG+H+h_{2}^{\prime\prime},\\
(ay+h_{1})H_{y}+(bz+u+h_{2})H_{z}+(bu+h_{3})H_{u}=bH+h_{3}^{\prime\prime}.
\end{array}\label{eq:oct7.2}
\end{equation}

\begin{itemize}
\item Case: $a\neq0,b\neq0$

Introducing the variables
\[
y=\overline{y}-\frac{h_{1}}{a},\ \ z=\overline{z}-\frac{h_{2}}{b}+\frac{h_{3}}{b^{2}},\ \ u=\overline{u}-\frac{h_{3}}{b},
\]
\[
F=\overline{F}-\frac{h_{1}^{\prime\prime}}{a},\ \ G=\overline{G}+\frac{h_{3}^{\prime\prime}}{b^{2}}-\frac{h_{2}^{\prime\prime}}{b},\ \ H=\overline{H}-\frac{h_{3}^{\prime\prime}}{b}
\]
 equations (\ref{eq:oct7.2}) become
\[
\begin{array}{l}
a\overline{y}\overline{F}_{\overline{y}}+(b\overline{z}+\overline{u})\overline{F}_{\overline{z}}+b\overline{u}\overline{F}_{\overline{u}}=a\overline{F},\\
a\overline{y}\overline{G}_{\overline{y}}+(b\overline{z}+\overline{u})\overline{G}_{\overline{z}}+b\overline{u}\overline{G}_{\overline{u}}=b\overline{G}+\overline{H},\\
a\overline{y}\overline{H}_{\overline{y}}+(b\overline{z}+\overline{u})\overline{H}_{\overline{z}}+b\overline{u}\overline{H}_{\overline{u}}=b\overline{H}.
\end{array}
\]
 In the variables
\[
\overline{y}=e^{as},\ \ \ \overline{z}=e^{bs}(sw+v),\ \ \ \overline{u}=e^{bs}w
\]
 these equations are
\[
\overline{F}_{s}=a\overline{F},\quad\overline{G}_{s}=b\overline{G}+c\overline{H},\quad\overline{H}_{s}=b\overline{H}
\]
 with general solution
\[
\overline{F}(x,y,z,u)=e^{as}f(x,v,w),
\]
\[
\overline{G}(x,y,z,u)=e^{bs}(h(x,v,w)s+g(x,v,w)),
\]
\[
\overline{H}(x,y,z,u)=e^{bs}h(x,v,w)
\]
 where $f(x,v,w),g(x,v,w)$ and $h(x,v,w)$ are arbitrary functions.

\end{itemize}

\begin{itemize}
\item Case: $a\neq0,b=0$

In this case equations (\ref{eq:oct7.2}) become
\begin{equation}
\begin{array}{l}
(ay+h_{1})F_{y}+(u+h_{2})F_{z}+h_{3}F_{u}=aF+h_{1}^{\prime\prime},\\
(ay+h_{1})G_{y}+(u+h_{2})G_{z}+h_{3}G_{u}=H+h_{2}^{\prime\prime},\\
(ay+h_{1})H_{y}+(u+h_{2})H_{z}+h_{3}H_{u}=h_{3}^{\prime\prime}.
\end{array}\label{cl:11}
\end{equation}
 Introducing the variables
\[
y=\overline{y}-\frac{h_{1}}{a},\ \ z=\overline{z},\ \ u=\overline{u}-h_{2},
\]
\[
F=\overline{F}-\frac{h_{1}^{\prime\prime}}{a},\ \ G=\overline{G},\ \ H=\overline{H}-h_{2}^{\prime\prime}
\]
 equations (\ref{cl:11}) become
\begin{equation}
\begin{array}{l}
a\overline{y}\overline{F}_{\overline{y}}+\overline{u}\overline{F}_{\overline{z}}+h_{3}\overline{F}_{\overline{u}}=a\overline{F},\\
a\overline{y}\overline{G}_{\overline{y}}+\overline{u}\overline{G}_{\overline{z}}+h_{3}\overline{G}_{\overline{u}}=\overline{H},\\
a\overline{y}\overline{H}_{\overline{y}}+\overline{u}\overline{H}_{\overline{z}}+h_{3}\overline{H}_{\overline{u}}=h_{3}^{\prime\prime}.
\end{array}\label{cl:16}
\end{equation}

\begin{itemize}
\item If $h_{3}\neq0$, then in the variables
\[
\overline{y}=e^{as},\ \ \ \overline{z}=\frac{h_{3}s^{2}}{2}+ws+v,\ \ \ \overline{u}=h_{3}s+w
\]
 equations (\ref{cl:16}) become
\[
\overline{F}_{s}=a\overline{F},\quad\overline{G}_{s}=\overline{H},\quad\overline{H}_{s}=h_{3}^{\prime\prime}
\]
 with general solution
\[
\overline{F}(x,y,z,u)=e^{as}f(x,v,w),
\]
\[
\overline{G}(x,y,z,u)=\frac{h_{3}s^{2}}{2}+h(x,v,w)s+g(x,v,w),
\]
\[
\overline{H}(x,y,z,u)=h_{3}^{\prime\prime}s+h(x,v,w).
\]

\end{itemize}

\begin{itemize}
\item If $h_{3}=0$ then in the variables
\[
\overline{y}=e^{as},\ \ \ \overline{z}=ws+v,\ \ \ \overline{u}=w
\]
 equations (\ref{cl:16}) become
\[
\overline{F}_{s}=a\overline{F},\quad\overline{G}_{s}=\overline{H},\quad\overline{H}_{s}=0
\]
 and their general solution is
\[
\overline{F}(x,y,z,u)=e^{as}f(x,v,w),
\]
\[
\overline{G}(x,y,z,u)=h(x,v,w)s+g(x,v,w),
\]
\[
\overline{H}(x,y,z,u)=h(x,v,w).
\]

\end{itemize}
\end{itemize}

\begin{itemize}
\item Case: $a=0,b\neq0$

In this case equations (\ref{eq:oct7.2}) become
\begin{equation}
\begin{array}{l}
h_{1}F_{y}+(bz+u+h_{2})F_{z}+(bu+h_{3})F_{u}=h_{1}^{\prime\prime},\\
h_{1}G_{y}+(bz+u+h_{2})G_{z}+(bu+h_{3})G_{u}=bG+H+h_{2}^{\prime\prime},\\
h_{1}H_{y}+(bz+u+h_{2})H_{z}+(bu+h_{3})H_{u}=bH+h_{3}^{\prime\prime}.
\end{array}\label{cl:13}
\end{equation}
 Introducing the variables
\[
y=\overline{y},\ \ z=\overline{z}-\frac{h_{2}}{b}+\frac{h_{3}}{b^{2}},\ \ u=\overline{u}-\frac{h_{3}}{b},
\]
\[
F=\overline{F},\ \ G=\overline{G}-\frac{h_{2}^{\prime\prime}}{b}+\frac{h_{3}^{\prime\prime}}{b^{2}},\ \ H=\overline{H}-\frac{h_{3}^{\prime\prime}}{b}
\]
 equations (\ref{cl:13}) become
\begin{equation}
\begin{array}{l}
h_{1}\overline{F}_{\overline{y}}+(b\overline{z}+\overline{u})\overline{F}_{\overline{z}}+b\overline{u}\overline{F}_{\overline{u}}=h_{1}^{\prime\prime},\\
h_{1}\overline{G}_{\overline{y}}+(b\overline{z}+\overline{u})\overline{G}_{\overline{z}}+b\overline{u}\overline{G}_{\overline{u}}=b\overline{G}+\overline{H},\\
h_{1}\overline{H}_{\overline{y}}+(b\overline{z}+\overline{u})\overline{H}_{\overline{z}}+b\overline{u}\overline{H}_{\overline{u}}=b\overline{H}.
\end{array}\label{cl:17}
\end{equation}

\begin{itemize}
\item If $h_{1}\neq0$, then in the variables
\[
\overline{y}=h_{1}s,\ \ \ \overline{z}=e^{bs}(ws+v),\ \ \ \overline{u}=e^{bs}w
\]
 equations (\ref{cl:17}) become
\[
\overline{F}_{s}=h_{1}^{\prime\prime},\quad\overline{G}_{s}=b\overline{G}+\overline{H},\quad\overline{H}_{s}=b\overline{H}.
\]
 The general solution of the last set of equations is
\[
\overline{F}(x,y,z,u)=h_{1}^{\prime\prime}s+f(x,v,w),
\]
\[
\overline{G}(x,y,z,u)=e^{bs}(h(x,v,w)s+g(x,v,w)),
\]
\[
\overline{H}(x,y,z,u)=e^{bs}h(x,v,w).
\]

\end{itemize}

\begin{itemize}
\item If $h_{1}=0$, then in the variables
\[
\overline{y}=w,\ \ \ \overline{z}=e^{bs}(w+v),\ \ \ \overline{u}=e^{bs}w
\]
 equations (\ref{cl:17}) become
\[
\overline{F}_{s}=0,\quad\overline{G}_{s}=b\overline{G}+\overline{H},\quad\overline{H}_{s}=b\overline{H}
\]
 and their general solution is
\[
\overline{F}(x,y,z,u)=f(x,v,w),
\]
\[
\overline{G}(x,y,z,u)=e^{bs}(h(x,v,w)s+g(x,v,w)),
\]
\[
\overline{H}(x,y,z,u)=e^{bs}h(x,v,w).
\]

\end{itemize}
\end{itemize}

\begin{itemize}
\item Case: $a=0,b=0$

In this case equations (\ref{eq:oct7.2}) become
\begin{equation}
\begin{array}{l}
h_{1}F_{y}+(u+h_{2})F_{z}+h_{3}F_{u}=h_{1}^{\prime\prime},\\
h_{1}G_{y}+(u+h_{2})G_{z}+h_{3}G_{u}=H+h_{2}^{\prime\prime},\\
h_{1}H_{y}+(u+h_{2})H_{z}+h_{3}H_{u}=h_{3}^{\prime\prime}.
\end{array}\label{cl:14}
\end{equation}
 Introducing the variables
\[
y=\overline{y},\ \ z=\overline{z},\ \ u=\overline{u}-h_{2},
\]
\[
F=\overline{F},\ \ G=\overline{G},\ \ H=\overline{H}-h_{2}^{\prime\prime}
\]
 equations (\ref{cl:14}) become
\begin{equation}
\begin{array}{l}
h_{1}\overline{F}_{\overline{y}}+\overline{u}\overline{F}_{\overline{z}}+h_{3}\overline{F}_{\overline{u}}=h_{1}^{\prime\prime},\\
h_{1}\overline{G}_{\overline{y}}+\overline{u}\overline{G}_{\overline{z}}+h_{3}\overline{G}_{\overline{u}}=\overline{H},\\
h_{1}\overline{H}_{\overline{y}}+\overline{u}\overline{H}_{\overline{z}}+h_{3}\overline{H}_{\overline{u}}=h_{3}^{\prime\prime}.
\end{array}\label{cl:15}
\end{equation}

\begin{itemize}
\item If $h_{1}\neq0,h_{3}\neq0$, then in the variables
\[
\overline{y}=h_{1}s,\ \ \ \overline{z}=\frac{h_{3}s^{2}}{2}+ws+v,\ \ \ \overline{u}=h_{3}s+w
\]
 equations (\ref{cl:15}) become
\[
\overline{F}_{s}=h_{1}^{\prime\prime},\quad\overline{G}_{s}=\overline{H},\quad\overline{H}_{s}=h_{3}^{\prime\prime}.
\]
 The general solution of the last set of equations is
\[
\overline{F}(x,y,z,u)=h_{1}^{\prime\prime}s+f(x,v,w),
\]
\[
\overline{G}(x,y,z,u)=\frac{h_{3}^{\prime\prime}s^{2}}{2}+h(x,v,w)s+g(x,v,w),
\]
\[
\overline{H}(x,y,z,u)=h_{3}^{\prime\prime}s+h(x,v,w).
\]

\end{itemize}

\begin{itemize}
\item If $h_{1}\neq0,h_{3}=0$, then in the variables
\[
\overline{y}=h_{1}s,\ \ \ \overline{z}=ws+v,\ \ \ \overline{u}=w
\]
 equations (\ref{cl:15}) become
\[
\overline{F}_{s}=h_{1}^{\prime\prime},\quad\overline{G}_{s}=\overline{H},\quad\overline{H}_{s}=0.
\]
 The general solution of the last set of equations is
\[
\overline{F}(x,y,z,u)=h_{1}^{\prime\prime}s+f(x,v,w),
\]
\[
\overline{G}(x,y,z,u)=h(x,v,w)s+g(x,v,w),
\]
\[
\overline{H}(x,y,z,u)=h(x,v,w).
\]

\end{itemize}

\begin{itemize}
\item If $h_{1}=0,h_{3}\neq0$, then in the variables
\[
\overline{y}=w,\ \ \ \overline{z}=\frac{h_{3}s^{2}}{2}+ws+v,\ \ \ \overline{u}=h_{3}s+w
\]
 equations (\ref{cl:15}) become
\[
\overline{F}_{s}=0,\quad\overline{G}_{s}=\overline{H},\quad\overline{H}_{s}=h_{3}^{\prime\prime}
\]
 with general solution
\[
\overline{F}(x,y,z,u)=f(x,v,w),
\]
\[
\overline{G}(x,y,z,u)=\frac{h_{3}^{\prime\prime}s^{2}}{2}+h(x,v,w)s+g(x,v,w),
\]
\[
\overline{H}(x,y,z,u)=h_{3}^{\prime\prime}s+h(x,v,w).
\]

\end{itemize}

\begin{itemize}
\item If $h_{1}=0,h_{3}=0$, then in the variables
\[
\overline{y}=v,\ \ \ \overline{z}=ws+v,\ \ \ \overline{u}=w
\]
 equations (\ref{cl:15}) become
\[
\overline{F}_{s}=0,\quad\overline{G}_{s}=\overline{H},\quad\overline{H}_{s}=0
\]
 and their general solution is
\[
\overline{F}(x,y,z,u)=f(x,v,w),
\]
\[
\overline{G}(x,y,z,u)=h(x,v,w)s+g(x,v,w),
\]
\[
\overline{H}(x,y,z,u)=h(x,v,w).
\]

\end{itemize}
\end{itemize}

\subsubsection{Case $A=J_{4}$}

In this case equations (\ref{eq:oct5.5}) become
\begin{equation}
\begin{array}{l}
(ay+z+h_{1})F_{y}+(az+u+h_{2})F_{z}+(au+h_{3})F_{u}=aF+G+h_{1}^{\prime\prime},\\
(ay+z+h_{1})G_{y}+(az+u+h_{2})G_{z}+(au+h_{3})G_{u}=aG+H+h_{2}^{\prime\prime},\\
(ay+z+h_{1})H_{y}+(az+u+h_{2})H_{z}+(au+h_{3})H_{u}=aH+h_{3}^{\prime\prime}.
\end{array}\label{cl:18}
\end{equation}

\begin{itemize}
\item Case: $a\neq0$

Introducing the variables
\[
y=\overline{y}-\frac{h_{1}}{a}+\frac{h_{2}}{a^{2}}-\frac{h_{3}}{a^{3}},\ \ z=\overline{z}-\frac{h_{2}}{a}+\frac{h_{3}}{a^{2}},\ \ u=\overline{u}-\frac{h_{3}}{a},
\]
\[
F=\overline{F}-\frac{h_{1}^{\prime\prime}}{a}+\frac{h_{2}^{\prime\prime}}{a^{2}}-\frac{h_{3}^{\prime\prime}}{a^{3}},\ \ G=\overline{G}-\frac{h_{2}^{\prime\prime}}{a}+\frac{h_{3}^{\prime\prime}}{a^{2}},\ \ H=\overline{H}-\frac{h_{3}^{\prime\prime}}{a}
\]
 equations (\ref{cl:18}) become
\[
\begin{array}{l}
(a\overline{y}+\overline{z})\overline{F}_{\overline{y}}+(a\overline{z}+\overline{u})\overline{F}_{\overline{z}}+a\overline{u}\overline{F}_{\overline{u}}=a\overline{F}+\overline{G},\\
(a\overline{y}+\overline{z})\overline{G}_{\overline{y}}+(a\overline{z}+\overline{u})\overline{G}_{\overline{z}}+a\overline{u}\overline{G}_{\overline{u}}=a\overline{G}+\overline{H},\\
(a\overline{y}+\overline{z})\overline{H}_{\overline{y}}+(a\overline{z}+\overline{u})\overline{H}_{\overline{z}}+a\overline{u}\overline{H}_{\overline{u}}=a\overline{H}.
\end{array}
\]
 In the variables
\[
\overline{y}=e^{as}(\frac{ws^{2}}{2}+v),\ \ \ \overline{z}=e^{as}(ws+v),\ \ \ \overline{u}=e^{as}w
\]
 these equations are
\[
\overline{F}_{s}=a\overline{F}+\overline{G},\quad\overline{G}_{s}=a\overline{G}+\overline{H},\quad\overline{H}_{s}=a\overline{H}.
\]
 The general solution of the last set of equations is
\[
\overline{F}(x,y,z,u)=e^{as}(h(x,v,w)\frac{s^{2}}{2}+g(x,v,w)s+f(x,v,w)),
\]
\[
\overline{G}(x,y,z,u)=e^{as}(h(x,v,w)s+g(x,v,w)),
\]
\[
\overline{H}(x,y,z,u)=e^{as}h(x,v,w)
\]
 where $f(x,v,w),g(x,v,w)$ and $h(x,v,w)$ are arbitrary functions.

\end{itemize}

\begin{itemize}
\item Case: $a=0$

In this case (\ref{cl:18}) become
\begin{equation}
\begin{array}{l}
(z+h_{1})F_{y}+(u+h_{2})F_{z}+h_{3}F_{u}=G+h_{1}^{\prime\prime},\\
(z+h_{1})G_{y}+(u+h_{2})G_{z}+h_{3}G_{u}=H+h_{2}^{\prime\prime},\\
(z+h_{1})H_{y}+(u+h_{2})H_{z}+h_{3}H_{u}=h_{3}^{\prime\prime}.
\end{array}\label{cl:19}
\end{equation}
 Introducing the variables
\[
y=\overline{y},\ \ z=\overline{z}-h_{1},\ \ u=\overline{u}-h_{2},
\]
\[
F=\overline{F},\ \ G=\overline{G}-h_{1}^{\prime\prime},\ \ H=\overline{H}-h_{2}^{\prime\prime}
\]
 equations (\ref{cl:19}) become
\begin{equation}
\begin{array}{l}
\overline{z}\overline{F}_{\overline{y}}+\overline{u}\overline{F}_{\overline{z}}+h_{3}\overline{F}_{\overline{u}}=\overline{G},\\
\overline{z}\overline{G}_{\overline{y}}+\overline{u}\overline{G}_{\overline{z}}+h_{3}\overline{G}_{\overline{u}}=\overline{H},\\
\overline{z}\overline{H}_{\overline{y}}+\overline{u}\overline{H}_{\overline{z}}+h_{3}\overline{H}_{\overline{u}}=h_{3}^{\prime\prime}.
\end{array}\label{cl:20}
\end{equation}

\begin{itemize}
\item If $h_{3}\neq0$, then in the variables
\[
\overline{y}=\frac{h_{3}s^{3}}{6}+vs+w,\ \ \ \overline{z}=\frac{h_{3}s^{2}}{2}+v,\ \ \ \overline{u}=h_{3}s
\]
 equations (\ref{cl:20}) become
\[
\overline{F}_{s}=\overline{G},\quad\overline{G}_{s}=\overline{H},\quad\overline{H}_{s}=h_{3}^{\prime\prime}
\]
 with general solution
\[
\overline{F}(x,y,z,u)=\frac{h_{3}s^{3}}{6}+h(x,v,w)\frac{s^{2}}{2}+g(x,v,w)s+f(x,v,w),
\]
\[
\overline{G}(x,y,z,u)=\frac{h_{3}s^{2}}{2}+h(x,v,w)s+g(x,v,w),
\]
\[
\overline{H}(x,y,z,u)=h_{3}s+h(x,v,w).
\]

\end{itemize}

\begin{itemize}
\item If $h_{3}=0$, then in the variables
\[
\overline{y}=\frac{ws^{2}}{2}+vs,\ \ \ \overline{z}=ws+v,\ \ \ \overline{u}=w
\]
 equations (\ref{cl:20}) become
\[
\overline{F}_{s}=\overline{G},\quad\overline{G}_{s}=\overline{H},\quad\overline{H}_{s}=0
\]
 with general solution
\[
\overline{F}(x,y,z,u)=h(x,v,w)\frac{s^{2}}{2}+g(x,v,w)s+f(x,v,w),
\]
\[
\overline{G}(x,y,z,u)=h(x,v,w)s+g(x,v,w),
\]
\[
\overline{H}(x,y,z,u)=h(x,v,w).
\]

\end{itemize}
\end{itemize}

\section{Systems of linear equations}

Linear second-order ordinary differential equations have the following
form,
\begin{equation}
\mathbf{y}^{\prime\prime}=A(x)\mathbf{y}^{\prime}+B(x)\mathbf{y}+f(x),\label{eq:mainlinear}
\end{equation}
 where $A(x)$ and $B(x)$ are matrices, and $f(x)$ is a vector.
Using a particular solution $\mathbf{y}_{p}(x)$ and the change
\[
\mathbf{y}=\widetilde{\mathbf{y}}+\mathbf{y}_{p},
\]
 we can, without loss of generality assume that $f(x)=0$. The matrix
$A(x)$ or $B(x)$ can also be assumed to be zero if we use the change,
$\mathbf{y}=C(x)\widetilde{\mathbf{y}}$, where $C=C(x)$ is a nonsingular
matrix. In the present paper the matrix $A(x)$ is reduced to zero.
In this case the linear equations (\ref{eq:main}) are linear functions
of $y,z$ and $u$:
\begin{equation}
\begin{array}{l}
F(x,y,z,u)=c_{11}(x)y+c_{12}(x)z+c_{13}(x)u,\\
G(x,y,z,u)=c_{21}(x)y+c_{22}(x)z+c_{23}(x)u,\\
H(x,y,z,u)=c_{31}(x)y+c_{32}(x)z+c_{33}(x)u.
\end{array}\label{eq:oct3.10-1}
\end{equation}
 In matrix form one can write as
\begin{equation}
\mathbf{F}(x,\mathbf{y})=C(x)\mathbf{y}\label{eq:oct2_2013.2}
\end{equation}

Any linear system of equations admits the following generators
\begin{equation}
y\partial_{y}+z\partial_{z}+u\partial_{u},\label{eq:oct7.10}
\end{equation}
\begin{equation}
\zeta_{1}(x)\partial_{y}+\zeta_{2}(x)\partial_{z}+\zeta_{3}(x)\partial_{u}\label{eq:oct7.10-1}
\end{equation}
 where $\zeta_{1}(x),\zeta_{2}(x)$ and $\zeta_{3}(x)$ are solutions
of the equations:
\[
\zeta_{1}^{\prime\prime}=c_{11}\zeta_{1}+c_{12}\zeta_{2}+c_{13}\zeta_{3},
\]
\[
\zeta_{2}^{\prime\prime}=c_{21}\zeta_{1}+c_{22}\zeta_{2}+c_{23}\zeta_{3},
\]
\[
\zeta_{3}^{\prime\prime}=c_{31}\zeta_{1}+c_{32}\zeta_{2}+c_{33}\zeta_{3}.
\]
 For the classification problem one needs to study equations which
admit generators different from (\ref{eq:oct7.10}) and (\ref{eq:oct7.10-1}).

\subsection{Equivalence transformations}

Similar to systems of two second-order ordinary differential equations
calculations show that the equivalence transformations are defined
by the following two types of transformations. The first type corresponds
to the linear change of the dependent variables $\widetilde{\mathbf{y}}=P\mathbf{y}$
with a constant nonsingular matrix $P$. The second type is
\[
\widetilde{x}=\varphi(x),\ \widetilde{y}=y\psi(x),\ \widetilde{z}=z\psi(x),\ \widetilde{u}=u\psi(x),
\]
 where the functions $\varphi(x)$ and $\psi(x)$ satisfy the condition
\[
\frac{\varphi^{\prime\prime}}{\varphi^{\prime}}=2\frac{\psi^{\prime}}{\psi}.
\]

{\bf Remark}. In \cite{bk:MkhizeMoyoMeleshko[2014]} it is shown that using the equivalence transformations of the presented above form, one can reduce any system of homogeneous linear second-order ordinary differential equations to a system with a matrix $C$ such that\footnote{For a system with two equations this property was obtained in \cite{bk:WafoMahomed[2000]}.}
 $trac\ C=0$. In this case one can obtain that the coefficient $\xi(x)$ of the admitted generator related with the independent variable $x$  has the form
 \[
 \xi=k_2x^2+k_1x+k_0,
 \]
 where $k_i$, ($i=0,1,2$) are constant. This means that the admitted generators are defined up to constants. This property allows one to apply an algebraic approach for group classification. In algebraic approach the constants are defined from the algebraic properties of admitted Lie algebras. Substituting the constants into determining equations one obtains systems of ordinary differential equations for the entries of the matrix $C(x)$. Solving this system of equations one finds forms of the entries.
 For systems of two linear second-order ordinary
differential equations it shown in \cite{bk:MkhizeMoyoMeleshko[2014]} that these systems of ordinary differential equations for the entries can be solved. However, the complete study of the
problem of group classification using the algebraic approach led to
the study of more cases than those found in \cite{bk:MoyoMeleshkoOguis[2013]},
where Ovsiannikov's approach was applied.

\subsection{Case $\xi\neq0$}

Using the obtained general forms of equations admitting an infinitesimal
generator with $\xi\neq0$ linear systems of equations (\ref{cl:02}),
(\ref{cl:04}), (\ref{cl:06}) and (\ref{cl:08}) have the following
form:
\begin{equation}
\begin{array}{l}
f(s,v,w)=\alpha_{11}s+\alpha_{12}v+\alpha_{13}w,\\
g(s,v,w)=\alpha_{21}s+\alpha_{22}v+\alpha_{23}w,\\
h(s,v,w)=\alpha_{31}s+\alpha_{32}v+\alpha_{33}w.
\end{array}\label{eq:oct3.10}
\end{equation}

\subsubsection{Case $A=J_{1}$}

Substituting (\ref{eq:oct3.10}) into (\ref{eq:oct3.10-1}), the equations
involving $F,G$ and $H$ become

\[
\begin{array}{l}
F=\alpha_{11}y+e^{\alpha x}\alpha_{12}z+e^{\beta x}\alpha_{13}u,\\
G=e^{-\alpha x}\alpha_{21}y+\alpha_{22}z+e^{-(\alpha-\beta)x}\alpha_{23}u,\\
H=e^{-\beta x}\alpha_{31}y+e^{(\alpha-\beta)x}\alpha_{32}z+\alpha_{33}u
\end{array}
\]
 where $\alpha=a-b$ and $\beta=a-d$.

Without loss of generality one can assume that $\alpha_{12}=1$. We
can assume that $\alpha_{12}\neq0$ because for $\alpha_{12}=0$ we
have $\alpha_{13}\neq0$ otherwise the system is degenerate. Applying
the scaling of $y$ we conclude that $\alpha_{12}=1$. Thus the studied
linear system of equations becomes
\begin{equation}
\begin{array}{l}
F=\alpha_{11}y+e^{\alpha x}z+e^{\beta x}\alpha_{13}u,\\
G=e^{-\alpha x}\alpha_{21}y+\alpha_{22}z+e^{-(\alpha-\beta)x}\alpha_{23}u,\\
H=e^{-\beta x}\alpha_{31}y+e^{(\alpha-\beta)x}\alpha_{32}z+\alpha_{33}u.
\end{array}\label{eq:sys_1}
\end{equation}

\subsubsection{Case $A=J_{2}$}

Substituting (\ref{eq:oct3.10}) into the functions $F,G$ and $H$
one finds
\begin{eqnarray*}
F & = & \alpha_{11}y+e^{\left(a-b\right)x}\left(\alpha_{12}\cos\left(cx\right)+\alpha_{13}\sin\left(cx\right)\right)z+e^{\left(a-b\right)x}(\alpha_{13}\cos\left(cx\right)\\
 &  & -\alpha_{12}\sin\left(cx\right))u,\\
G & = & e^{-\left(a-b\right)x}(\alpha_{21}\cos\left(cx\right)+\alpha_{31}\sin\left(cx\right))y+(\alpha_{22}\cos^{2}\left(cx\right)\\
 &  & +(\alpha_{23}+\alpha_{32})\cos\left(cx\right)\sin\left(cx\right)+\alpha_{33}\sin^{2}\left(cx\right))z\\
 &  & +(\alpha_{23}\cos^{2}\left(cx\right)+(\alpha_{33}-\alpha_{22})\cos\left(cx\right)\sin\left(cx\right)-\alpha_{32}\sin^{2}\left(cx\right))u,\\
\\
H & = & e^{-\left(a-b\right)x}(\alpha_{31}\cos\left(cx\right)-\alpha_{21}\sin\left(cx\right))y+(\alpha_{32}\cos^{2}\left(cx\right)\\
 &  & +(\alpha_{33}-\alpha_{22})\cos\left(cx\right)\sin\left(cx\right)-\alpha_{23}\sin^{2}\left(cx\right))z\\
 &  & +(\alpha_{33}\cos^{2}\left(cx\right)-(\alpha_{32}+\alpha_{23})\cos\left(cx\right)\sin\left(cx\right)+\alpha_{22}\sin^{2}\left(cx\right))u.\\
\end{eqnarray*}
 Using trigonometry formulae and introducing the constants $\beta$,
$\gamma$, $c_{1}$ and $c_{2}$:
\[
\left(\begin{array}{cc}
\alpha_{22} & \alpha_{23}\\
\alpha_{32} & \alpha_{33}
\end{array}\right)=\left(\begin{array}{cc}
\beta+c_{2} & \gamma-c_{1}\\
\gamma+c_{1} & -\beta+c_{2}
\end{array}\right),
\]
 one can rewrite functions $F$, $G$ and $H$ in the form
\begin{eqnarray*}
F & = & \alpha_{11}y+e^{\alpha x}\left(\cos\left(cx\right)\alpha_{12}+\sin\left(cx\right)\alpha_{13}\right)z+e^{\alpha x}(\cos\left(cx\right)\alpha_{13}\\
 &  & -\sin\left(cx\right)\alpha_{12})u,\\
G & = & e^{-\alpha x}\left(\cos\left(cx\right)\alpha_{21}+\sin\left(cx\right)\alpha_{31}\right)y+(\cos\left(2cx\right)\beta\\
 &  & +\sin\left(2cx\right)\gamma+c_{2})z+\left(\cos\left(2cx\right)\gamma-\sin\left(2cx\right)\beta-c_{1}\right)u,\\
H & = & e^{-\alpha x}\left(\cos\left(cx\right)\alpha_{31}-\sin\left(cx\right)\alpha_{21}\right)y+(\cos\left(2cx\right)\gamma\\
 &  & -\sin\left(2cx\right)\beta+c_{1})z-\left(\cos\left(2cx\right)\beta+\sin\left(2cx\right)\gamma-c_{2}\right)u,
\end{eqnarray*}
 where $\alpha=a-b$.

Applying rotation of the dependent variables $z$ and $u$, and their
dilation, we can assume that $\alpha_{12}=1$ and $\alpha_{13}=0$.
Hence the studied linear system of equations becomes
\begin{equation}
\begin{array}{l}
F=\alpha_{11}y+e^{\alpha x}\cos\left(cx\right)z-e^{\alpha x}\sin\left(cx\right)u,\\
G=e^{-\alpha x}\left(\cos\left(cx\right)\alpha_{21}+\sin\left(cx\right)\alpha_{31}\right)y+(\cos\left(2cx\right)\beta\\
\;\;\;\;\;\;+\sin\left(2cx\right)\gamma+c_{2})z+\left(\cos\left(2cx\right)\gamma-\sin\left(2cx\right)\beta-c_{1}\right)u,\\
H=e^{-\alpha x}\left(\cos\left(cx\right)\alpha_{31}-\sin\left(cx\right)\alpha_{21}\right)y+(\cos\left(2cx\right)\gamma\\
\;\;\;\;\;\;-\sin\left(2cx\right)\beta+c_{1})z-\left(\cos\left(2cx\right)\beta+\sin\left(2cx\right)\gamma-c_{2}\right)u.
\end{array}\label{eq:sys_2}
\end{equation}

\subsubsection{Case $A=J_{3}$}

Substituting (\ref{eq:oct3.10}) into the functions $F,G$ and $H$
one finds
\begin{equation}
\begin{array}{l}
F=\alpha_{11}y+e^{\alpha x}\alpha_{12}z+e^{\alpha x}\left(-\alpha_{12}x+\alpha_{13}\right)u,\\
G=e^{-\alpha x}\left(\alpha_{21}+\alpha_{31}x\right)y+\left(\alpha_{22}+\alpha_{32}x\right)z\;\\
\;\;\;\;\;\;+\left(\alpha_{23}+\left(\alpha_{33}-\alpha_{22}\right)x-\alpha_{32}x^{2}\right)u,\\
H=e^{-\alpha x}\alpha_{31}y+\alpha_{32}z+\left(\alpha_{33}-\alpha_{32}x\right)u
\end{array}\label{eq:sys_3}
\end{equation}
 where $\alpha=a-b$.

\subsubsection{Case $A=J_{4}$}

Substituting (\ref{eq:oct3.10}) into the functions $F,G$ and $H$
one finds
\begin{equation}
\begin{array}{l}
F=\left(\lambda+\beta x+\frac{1}{2}\gamma x^{2}\right)\left(y-xz+\frac{1}{2}x^{2}u\right),\\
G=e^{-\alpha x}\left(\beta+\gamma x\right)\left(y-xz+\frac{1}{2}x^{2}u\right),\\
H=e^{-\alpha x}\gamma\left(y-xz+\frac{1}{2}x^{2}u\right)
\end{array}\label{eq:sys_4}
\end{equation}
 where $\lambda=\alpha_{11}+\alpha_{12}+\alpha_{13},\,\,\,\beta=\alpha_{21}+\alpha_{22}+\alpha_{23},\,\,\,\gamma=\alpha_{31}+\alpha_{32}+\alpha_{33}.$

Since for $\gamma=0$ the linear system of equations is reduced to
the degenerate case $H=0$, one has to consider $\gamma\neq0$.

\subsection{Case $\xi=0$}

Substituting (\ref{eq:oct2_2013.2}) into (\ref{eq:oct2_2013.1})
and splitting it with respect to $y,z$ and $u$ one has
\begin{equation}
CA-AC=0.\label{eq:oct8.1}
\end{equation}
 A nontrivial admitted generator is of the form
\begin{equation}
X_{o}=(A\mathbf{y})\cdot\nabla.\label{eq:oct8.5}
\end{equation}

Equations (\ref{eq:oct8.1}) can be simplified by using the Jordan
form of the matrix $A$.

In particular, if $A=J_{1}$, then equations (\ref{eq:oct8.1}) become
\begin{equation}
\begin{array}{l}
c_{12}\left(b-a\right)=0,\quad c_{13}\left(d-a\right)=0,\quad c_{21}\left(a-b\right)=0,\\
c_{23}\left(b-d\right)=0,\quad c_{31}\left(a-d\right)=0,\quad c_{32}\left(b-d\right)=0.
\end{array}\label{eq:oct2_2013.5}
\end{equation}
 If $\left(a-b\right)^{2}+\left(a-d\right)^{2}+\left(b-d\right)^{2}=0$,
then we find that $a=b=d$ and the generator (\ref{eq:oct8.5}) is
also trivial. Hence we can assume that $a\neq b$. Using this assumption
we see from equations (\ref{eq:oct2_2013.5}) that
\[
c_{12}=0,\ \ c_{21}=0.
\]
 For a non-degenerate linear system one has to assume that $c_{13}\neq0$,
which leads to $a=d$ and
\[
c_{23}=0,\: c_{32}=0.
\]
 The linear system of equations with $c_{21}=0$ and $c_{23}=0$ is
degenerate.

Similar results are obtained for the for other cases of the matrix
$A$. This means that in contrast to a linear system with two equations
there is no linear systems with three dependent variables such that
all nontrivial admitted generators have the form (\ref{eq:oct8.5}).

\section{Solutions of the determining equations}

Solving the determining equations leads to the following solutions:
\begin{itemize}
\item For system (\ref{eq:sys_1}) there is the only nontrivial admitted
generator
\begin{equation}
\partial_{x}-\alpha z\partial_{z}-\beta u\partial_{u}.\label{eq:sys_1_gen}
\end{equation}

\item In the case of system (\ref{eq:sys_2}) there is the only nontrivial admitted
generator
\begin{equation}
\partial_{x}+\alpha y\partial_{y}+c\left(u\partial_{z}-z\partial_{u}\right).\label{eq:sys_2_gen}
\end{equation}

\item Similarly for system (\ref{eq:sys_3}) there is the only nontrivial admitted
generator
\begin{equation}
\partial_{x}+\alpha y\partial_{y}+u\partial_{z}.\label{eq:sys_3_gen}
\end{equation}

\item Finally for system (\ref{eq:sys_4}) there is the only nontrivial admitted
generator
\begin{equation}
\partial_{x}+z\partial_{y}+u\partial_{z}.\label{eq:sys_4_gen}
\end{equation}

\end{itemize}

\textbf{Theorem}.
All non-degenerate linear systems of three second-order
ordinary differential equations admitting a non-trivial generator
are equivalent to one of the cases:
\\
(a) system (\ref{eq:sys_1}) with the generator (\ref{eq:sys_1_gen}),
\\
(b) system (\ref{eq:sys_2}) with the generator (\ref{eq:sys_2_gen}),
\\
(c) system (\ref{eq:sys_3}) with the generator (\ref{eq:sys_3_gen}),
\\
(d) system (\ref{eq:sys_4}) with the generator (\ref{eq:sys_4_gen}).

\section{Conclusion}

The results of the current paper are obtained using Ovsiannikov's
approach \cite{bk:Ovsiannikov[2004]}, which involves simplifying
one generator and finding the associated functions. These functions
are then used to solve the determining equations. We found all forms
of nonlinear systems $\mathbf{y}^{\prime\prime}=\mathbf{F}(x,\mathbf{y})$
admitting at least one generator. Finding their forms is reduced to
solving a linear system of first-order ordinary differential equations
with constant coefficients. Methods for solving these systems are
well-known in the theory of ordinary differential equations; their
solution is given by studying purely algebraic properties of some
constant matrix $A$. Since the matrix $A$ is considered in the Jordan
form, these properties are simple in each particular case.

Any normal linear system of second-order ordinary differential equations
can be reduced to the form
\[
\mathbf{y}^{\prime\prime}=C\mathbf{y},
\]
 where $C(x)$ is a square matrix. These systems admit the set of
trivial generators $\mathbf{y}\cdot\nabla$, $\mathbf{h}(x)\cdot\nabla,$
where $\mathbf{h}^{\prime\prime}=C\mathbf{h}$. Complete group classification
of linear systems containing three equations and admitting a nontrivial
generator is given; it is shown that all such systems have one of
the four forms . We have also listed an explicit form of the non-trivial
admitted generators for each of the respective cases. The obtained
results are summarized in the theorem.

\section*{Acknowledgements}

This research was financially supported by Thailand Research Fund
under Grant no. MRG5580053.

\end{document}